\magnification=1200

\hsize=11.25cm    
\vsize=18cm       
\parindent=12pt   \parskip=5pt     

\hoffset=.5cm   
\voffset=.8cm   

\pretolerance=500 \tolerance=1000  \brokenpenalty=5000

\catcode`\@=11

\font\eightrm=cmr8         \font\eighti=cmmi8
\font\eightsy=cmsy8        \font\eightbf=cmbx8
\font\eighttt=cmtt8        \font\eightit=cmti8
\font\eightsl=cmsl8        \font\sixrm=cmr6
\font\sixi=cmmi6           \font\sixsy=cmsy6
\font\sixbf=cmbx6

\font\tengoth=eufm10 
\font\eightgoth=eufm8  
\font\sevengoth=eufm7      
\font\sixgoth=eufm6        \font\fivegoth=eufm5

\skewchar\eighti='177 \skewchar\sixi='177
\skewchar\eightsy='60 \skewchar\sixsy='60

\newfam\gothfam           \newfam\bboardfam

\def\tenpoint{
  \textfont0=\tenrm \scriptfont0=\sevenrm \scriptscriptfont0=\fiverm
  \def\rm{\fam\z@\tenrm}
  \textfont1=\teni  \scriptfont1=\seveni  \scriptscriptfont1=\fivei
  \def\oldstyle{\fam\@ne\teni}\let\old=\oldstyle
  \textfont2=\tensy \scriptfont2=\sevensy \scriptscriptfont2=\fivesy
  \textfont\gothfam=\tengoth \scriptfont\gothfam=\sevengoth
  \scriptscriptfont\gothfam=\fivegoth
  \def\goth{\fam\gothfam\tengoth}
  
  \textfont\itfam=\tenit
  \def\it{\fam\itfam\tenit}
  \textfont\slfam=\tensl
  \def\sl{\fam\slfam\tensl}
  \textfont\bffam=\tenbf \scriptfont\bffam=\sevenbf
  \scriptscriptfont\bffam=\fivebf
  \def\bf{\fam\bffam\tenbf}
  \textfont\ttfam=\tentt
  \def\tt{\fam\ttfam\tentt}
  \abovedisplayskip=12pt plus 3pt minus 9pt
  \belowdisplayskip=\abovedisplayskip
  \abovedisplayshortskip=0pt plus 3pt
  \belowdisplayshortskip=4pt plus 3pt 
  \smallskipamount=3pt plus 1pt minus 1pt
  \medskipamount=6pt plus 2pt minus 2pt
  \bigskipamount=12pt plus 4pt minus 4pt
  \normalbaselineskip=12pt
  \setbox\strutbox=\hbox{\vrule height8.5pt depth3.5pt width0pt}
  \let\bigf@nt=\tenrm       \let\smallf@nt=\sevenrm
  \normalbaselines\rm}

\def\eightpoint{
  \textfont0=\eightrm \scriptfont0=\sixrm \scriptscriptfont0=\fiverm
  \def\rm{\fam\z@\eightrm}
  \textfont1=\eighti  \scriptfont1=\sixi  \scriptscriptfont1=\fivei
  \def\oldstyle{\fam\@ne\eighti}\let\old=\oldstyle
  \textfont2=\eightsy \scriptfont2=\sixsy \scriptscriptfont2=\fivesy
  \textfont\gothfam=\eightgoth \scriptfont\gothfam=\sixgoth
  \scriptscriptfont\gothfam=\fivegoth
  \def\goth{\fam\gothfam\eightgoth}
  
  \textfont\itfam=\eightit
  \def\it{\fam\itfam\eightit}
  \textfont\slfam=\eightsl
  \def\sl{\fam\slfam\eightsl}
  \textfont\bffam=\eightbf \scriptfont\bffam=\sixbf
  \scriptscriptfont\bffam=\fivebf
  \def\bf{\fam\bffam\eightbf}
  \textfont\ttfam=\eighttt
  \def\tt{\fam\ttfam\eighttt}
  \abovedisplayskip=9pt plus 3pt minus 9pt
  \belowdisplayskip=\abovedisplayskip
  \abovedisplayshortskip=0pt plus 3pt
  \belowdisplayshortskip=3pt plus 3pt 
  \smallskipamount=2pt plus 1pt minus 1pt
  \medskipamount=4pt plus 2pt minus 1pt
  \bigskipamount=9pt plus 3pt minus 3pt
  \normalbaselineskip=9pt
  \setbox\strutbox=\hbox{\vrule height7pt depth2pt width0pt}
  \let\bigf@nt=\eightrm     \let\smallf@nt=\sixrm
  \normalbaselines\rm}

\tenpoint

\def\pc#1{\bigf@nt#1\smallf@nt}         \def\pd#1 {{\pc#1} }

\catcode`\;=\active
\def;{\relax\ifhmode\ifdim\lastskip>\z@\unskip\fi
\kern\fontdimen2  -1.2 \fontdimen3 \string;}

\catcode`\:=\active
\def:{\relax\ifhmode\ifdim\lastskip>\z@\unskip\fi\penalty\@M\ \fi\string:}

\catcode`\!=\active
\def!{\relax\ifhmode\ifdim\lastskip>\z@
\unskip\fi\kern\fontdimen2  -1.1 \fontdimen3 \string!}

\catcode`\?=\active
\def?{\relax\ifhmode\ifdim\lastskip>\z@
\unskip\fi\kern\fontdimen2  -1.1 \fontdimen3 \string?}

\frenchspacing

\def\raggedbottom{\topskip 10pt plus 36pt\r@ggedbottomtrue}

\def\pointir{\unskip . --- \ignorespaces}

\def\Medbreak{\vskip-\lastskip\medbreak}

\long\def\th#1 #2\enonce#3\endth{
   \Medbreak\noindent
   {\pc#1} {#2\unskip}\pointir{\it #3}\smallskip}

\def\decale#1{\smallbreak\hskip 28pt\llap{#1}\kern 5pt}
\def\decaledecale#1{\smallbreak\hskip 34pt\llap{#1}\kern 5pt}
\def\puce{\smallbreak\hskip 6pt{$\scriptstyle\bullet$}\kern 5pt}

\def\eqalign#1{\null\,\vcenter{\openup\jot\m@th\ialign{
\strut\hfil$\displaystyle{##}$&$\displaystyle{{}##}$\hfil
&&\quad\strut\hfil$\displaystyle{##}$&$\displaystyle{{}##}$\hfil
\crcr#1\crcr}}\,}

\catcode`\@=12

\showboxbreadth=-1  \showboxdepth=-1

\newcount\numerodesection \numerodesection=1
\def\section#1{\bigbreak
 {\bf\number\numerodesection.\ \ #1}\nobreak\medskip
 \advance\numerodesection by1}

\mathcode`A="7041 \mathcode`B="7042 \mathcode`C="7043 \mathcode`D="7044
\mathcode`E="7045 \mathcode`F="7046 \mathcode`G="7047 \mathcode`H="7048
\mathcode`I="7049 \mathcode`J="704A \mathcode`K="704B \mathcode`L="704C
\mathcode`M="704D \mathcode`N="704E \mathcode`O="704F \mathcode`P="7050
\mathcode`Q="7051 \mathcode`R="7052 \mathcode`S="7053 \mathcode`T="7054
\mathcode`U="7055 \mathcode`V="7056 \mathcode`W="7057 \mathcode`X="7058
\mathcode`Y="7059 \mathcode`Z="705A


\def\hfl#1#2#3{\smash{\mathop{\hbox to#3{\rightarrowfill}}\limits
^{\textstyle#1}_{\textstyle#2}}}

\def\Q{{\bf Q}}

\def\R{{\bf R}}

\def\Z{{\bf Z}}

\def\F{{\bf F}}

\def\Gal{\mathop{\rm Gal}\nolimits}

\def\to{\rightarrow}

\def\normressym(#1,#2)_#3{\displaystyle\left({#1,#2\over#3}\right)}

\newcount\refno 
\long\def\ref#1:#2<#3>{                                        
\global\advance\refno by1\par\noindent                              
\llap{[{\bf\number\refno}]\ }{#1} \pointir{\it #2} #3\goodbreak }

\def\citer#1(#2){[{\bf\number#1}\if#2\empty\relax\else,\ {#2}\fi]}

\newbox\bibbox
\setbox\bibbox\vbox{\bigbreak
\centerline{{\pc BIBLIOGRAPHY}}

\ref{\pc CAPUANO} (L) \& {\pc DEL \pc CORSO} (I):
Upper ramification jumps in abelian extensions of exponent~$p$,
<Preprint, arXiv\string:1407.2496.>
\newcount\capuanocorso \global\capuanocorso=\refno

\ref{\pc CASTA\~NEDA} (J) \& {\pc WU} (Q): 
The ramification group filtrations of certain function field
extensions,
<Pacific Journal of Mathematics 276.2 (2015), 309--320.>
\newcount\castawu \global\castawu=\refno

\ref{\pc DALAWAT} (C):
Local discriminants, kummerian extensions, and elliptic curves,
<J. Ramanujan Math.\ Soc.\ 25.1 (2010),
25--80. Cf.~arXiv\string:0711.3878.>      
\newcount\locdisc \global\locdisc=\refno

\ref{\pc DALAWAT} (C):
Further remarks on local discriminants, 
<J. Ramanujan Math.\ Soc.\ 25.4 (2010), 393--417. Cf.~arXiv\string:0909.2541.>    
\newcount\further \global\further=\refno

\ref{\pc DALAWAT} (C):
Final remarks on local discriminants, 
<J. Ramanujan Math.\ Soc.\ 25.4 (2010), 419--432. Cf.~arXiv\string:0912.2829.>    
\newcount\final \global\final=\refno

\ref{\pc SERRE} (J-P):
Corps locaux,
<Publications de l'Universit{\'e} de Nancago, No.~{\sevenrm VIII}, Hermann,
Paris, 1968, 245 pp.>
\newcount\corpslocaux \global\corpslocaux=\refno

} 

\centerline{\bf The ramification filtration in certain $p$-extensions} 
\bigskip\bigskip 
\centerline{Chandan Singh Dalawat} 
\centerline{Harish-Chandra Research Institute}
\centerline{Chhatnag Road, Jhunsi, Allahabad 211019, India} 
\centerline{\tt dalawat@gmail.com}

\bigskip\bigskip

{{\bf Abstract}.  We show that the recent result of Casta\~neda and Wu
about the ramification filtration in certain $p$-extensions of
function fields of prime characteristic $p$ is equally valid over
local fields of mixed characteristic $(0,p)$.  Apart from being
applicable to both equicharacteristic and mixed characteristic cases,
our method has the advantage of being purely local, purely
conceptual, more natural, and much shorter.

\footnote{}{{\it MSC2010~:} primary 11S15}
\footnote{}{{\it Keywords~:} Local fields,
    abelian extensions of exponent~$p$, ramification filtration, upper
ramification breaks}}

\medskip
\rightline{\it Il faut le faire.}
\rightline{--- Samuel Eilenberg}
\smallskip

{\bf 1.  Introduction}\pointir Let $p$ be a prime number and $K$ a
local field with perfect residue field of characteristic~$p$.  Let $G$
be a $p$-group which has the property that every subgroup of $G$ is
the intersection of the family of index-$p$ subgroups of $G$
containing it, and let $L$ be a totally ramified $G$-extension of $K$.
When $K$ has characteristic~$p$, Casta\~neda and
Wu \citer\castawu(Theorem~4.4) have recently established a
relationship between the possible exponents of the differents of
intermediate extensions $L|E|K$ which have degree~$p$ over $K$ and the
lower ramification breaks of $L|K$.  We show that the same
relationship holds even when $K$ has characteristic~$0$, and it is
more easily derived by using the ramification filtration in the upper
numbering and Herbrand's theorem.  We also remark that the hypothesis
on the group $G$ forces it to be commutative of exponent~$p$.

\medbreak

{\bf 2.  Certain $p$-groups}\pointir Let us first show that the
$p$-groups which were introduced in \citer\castawu(2-5) are the same
as $\F_p$-spaces~:

\th LEMMA 2.1
\enonce
If\/ $G$ is a $p$-group\/ (of order~$>1$) in which every subgroup
$G'\subset G$ is the intersection of the family of index-$p$ subgroups
of\/ $G$ containing\/ $G'$, then\/ $G$ is commutative of exponent~$p$.
\endth

{\it Proof~:\/} Taking $G'$ to be the trivial subgroup of $G$, we
conclude that the intersection of the family of {\it all\/} index-$p$
subgroups of $G$ is trivial.  Let $P$ be the maximal quotient of $G$
which is commutative of exponent~$p$ and let $N$ be the kernel of the
projection $G\to P$~; we have to show that $N$ is trivial, for which
it is enough to show that $N\subset H$ for every index-$p$ subgroup
$H\subset G$, as we have seen.  Recall that every index-$p$ subgroup
in a $p$-group is normal.  As $G/H$ is commutative of exponent~$p$,
there is a unique morphism $P\to G/H$ such that the projection $G\to
G/H$ factors as $G\to P\to G/H$ (by the maximality of $P$), and hence
$N\subset H$.  Therefore $N$ is trivial, $G=P$, and the lemma is proved. 

\medbreak

{\bf 3.  The upper ramification breaks}\pointir It follows that
$G$-extensions $L|K$ for which the group $G$ satisfies the hypothesis
of Lemma~2.1 are the same as abelian extensions of exponent~$p$.  We
next determine the upper ramification breaks of such extensions in
terms of the degree-$p$ extensions of $K$ contained in $L$,
necessarily cyclic over $K$.

\th PROPOSITION 3.1
\enonce
For every abelian extension\/ $L$ of\/ $K$ of exponent\/~$p$, the
upper ramification breaks\/ $t$ of\/ $L|K$ are precisely the same as the
possible ramification breaks of intermediate extensions\/ $L|E|K$
which are cyclic and of degree~$p$ over\/~$K$.  In particular, $t\in\Z$.
\endth

{\it Proof~:\/} We will use the compatibility of the ramification
filtration in the upper numbering $(G^u)_{u\in\R}$ on $G=\Gal(L|K)$
with the passage to the quotient (Herbrand's
theorem \citer\corpslocaux(Chapter~IV, Proposition~14)).

Explicitly, denoting by $G^{u+}$ the union of $G^{u+\varepsilon}$
($\varepsilon>0$), if there is an intermediate extension $L|E|K$ which
is (cyclic) of degree~$p$ over~$K$, and if the unique ramification
break of $E|K$ occurs at $t$ (necessarily an integer), then $C^t=C$
and $C^{t+}=\{1\}$, where $C=\Gal(E|K)$.  But by Herbrand's theorem,
$C^t=\overline{G^t}$ and $C^{t+}=\overline{G^{t+}}$, where
$\overline{H}$ denotes the image in $C$ of a subgroup $H\subset G$
under the projection $G\to C$.  Hence $G^{t+}\neq G^t$, and an upper
ramification break occurs for $L|K$ at $t$. (The hypothesis that $G$
be an $\F_p$-space has not been used yet.)

For the converse, suppose that an upper ramification break occurs at
$t$ for $L|K$, so that $G^{t+}\neq G^t$.  Since $G$ is an
$\F_p$-space, there is an index-$p$ subgroup $H\subset G$ such that
$G^{t+}\subset H$ but $G^t\not\subset H$, so that $L^H$ is a
degree-$p$ (cyclic) extension of $K$.  For $C=\Gal(L^H|K)=G/H$, by
Herbrand's theorem we have $C^t=\overline{G^t}=C$ because
$G^t\not\subset H$, and $C^{t+}=\overline{G^{t+}}=\{1\}$ because
$G^{t+}\subset H$, so the unique ramification break for $L^H|K$ occurs
at $t$ (which therefore has to be an integer).  This completes the
proof.

{\it Remark~3.2}.  It is known that the upper ramification breaks of
{\it any\/} abelian extension of $K$ are integers
(Hasse-Arf, \citer\corpslocaux(Chapter~V, Theorem~1)).
Proposition~3.1 implies and sharpens Hasse-Arf for abelian extensions
of exponent~$p$ by specifying which integers can occur.

{\it Remark~3.3}.  In view of Proposition~3.1, one would like to
determine the possibilities for the unique ramification break $t$ of a
{\it ramified\/} degree-$p$ cyclic extension $E|K$.  For this purpose,
denote by
$$
b^{(i)}=i+\lfloor{(i-1)/(p-1)}\rfloor\qquad (i>0)
$$
the sequence of positive integers which are prime~to~$p$.  If $K$ has
characteristic~$0$, let $e$ be the ramification index of $K|\Q_p$.  If
moreover $K$ does not contain a primitive $p$-th root $\zeta$ of~$1$,
the possibilities for $t$ are $b^{(1)},\ldots,b^{(e)}$ and each of
them does occur \citer\locdisc(Proposition~63).  If $K$ contains
$\zeta$, there is one further possibility, namely $pe_1=b^{(e)}+1$,
where $e_1$ is the ramification index of $K|\Q_p(\zeta)$ (so that
$e=(p-1)e_1$), and it does occur for $E=K(\!\root p\of\pi)$, where $\pi$
is a uniformiser of~$K$ \citer\locdisc(Corollary~62).  If $K$ has
characteristic~$p$, each $b^{(i)}$ is a possibility, and each of them
does occur \citer\further(Proposition~14).  

{\it Remark~3.4}.  The special case of Proposition~3.1 when $K$ is a
finite extension of $\Q_p$ can be found in the recent preprint of
Capuano and Del Corso \citer\capuanocorso(Proposition~7).
Incidentally, \citer\further() and \citer\final() derive the results
of \citer\capuanocorso() without using class field theory and also
treat the equicharacteristic case.

{\it Remark~3.5}.  If the unique ramification break of a ramified
degree-$p$ cyclic extension $E$ of $K$ occurs at $t$, then the
exponent of the different of $E|K$ is
$d=(1+t)(p-1)$ \citer\corpslocaux(Chapter~IV, Proposition~4), so $t$
can be recovered from $d$ and Proposition~3.1 could have been stated
in terms of the possible $d$ instead of the possible $t$ for such
$E\subset L$.

\medbreak

{\bf 4.  The lower ramification breaks}\pointir Let us finally show
how to recover the main result of Casta\~neda and
Wu \citer\castawu(Theorem~4.4) (and extend it to the mixed
characteristic case) from the foregoing.

So let $L|K$ be a totally ramified abelian extension of exponent~$p$,
and assume that the degree~$[L:K]$ is finite, so that the filtration
in the lower numbering is defined on $G=\Gal(L|K)$.  Let
$$
d_1<d_2<\cdots<d_n
$$
be the exponents of the differents of $E|K$, for the various
degree-$p$ extensions $E$ of $K$ in $L$.  We have to compute the {\it
lower\/} ramification breaks of $L|K$ in terms of these $d_i$ and the
indices $\hbox{{\bf(}}G^t:G^{t+}\hbox{{\bf)}}$ for the various {\it
upper\/} ramification breaks $t$ of $L|K$.

We have seen in Remark~3.5 that the unique ramification breaks of
intermediate extensions $L|E|K$ which are cyclic and of degree~$p$
over $K$ occur precisely at
$$
t_1<t_2<\cdots<t_n,
$$
where $d_i=(1+t_i)(p-1)$ for every $i\in[1,n]$.  Note that $t_1>0$
because $L|K$ is a totally ramified $p$-extension (of degree~$>1$).

It follows from Proposition~3.1 that the upper ramification breaks of
$L|K$ occur precisely at these $t_i$.  For each $i\in[1,n]$, let $f_i$
be the codimension of $G^{t_i+}$ in $G^{t_i}$, or equivalently
$p^{f_i}=\hbox{{\bf(}}G^{t_i}:G^{t_i+}\hbox{{\bf)}}$ (so that
$f_i>0$).  Note that $G^{t_i+}=G^{t_{i+1}}$ for every $i\in[1,n[$.
All this information is neatly summarised by the diagram
$$
\{1\}
=G^{t_n+}
\subset_{f_n} G^{t_n}
\subset_{f_{n-1}} G^{t_{n-1}}
\subset_{f_{n-2}}
\cdots
\subset_{f_2} G^{t_2}
\subset_{f_1} G^{t_1}
=G
$$
where the notation $V\subset_c W$ means that $V$ is a codimension-$c$
subspace of the $\F_p$-space $W$.

\th PROPOSITION 4.1
\enonce
The lower ramification breaks of the extension\/ $L|K$ occur precisely
at\/ $l_1=t_1$ and at\/
$l_i=l_{i-1}+(t_i-t_{i-1})p^{f_1+f_2+\cdots+f_{i-1}}$ for\/ $i\in[2,n]$.
\endth

{\it Proof~:\/}  We use the formula $l_i=\psi_{L|K}(t_i)$ for the
passage from the upper numbering to the lower
numbering \citer\corpslocaux(Chapter~IV, \S3).  Recall that
$$
\psi_{L|K}(v)=\int_0^v\!\!\hbox{\bf(}G^0:G^w\hbox{\bf)}\,dw. 
$$
The result follows upon tabulating the indices
$\hbox{{\bf(}}G^0:G^w\hbox{\bf)}$ for every real $w\in[0,+\infty[$~:
$$
\vbox{\halign{&\hfil$#$\hfil\quad\cr
\multispan7\hrulefill\cr
w\in&[0,t_1]&]t_1,t_2]&]t_2,t_3]&
\cdots&]t_{n-1},t_n]&]t_n,+\infty[\cr
\noalign{\vskip-5pt}
\multispan7\hrulefill\cr
\hbox{\bf(}G^0:G^w\hbox{\bf)}=
&1&p^{f_1}&p^{f_1+f_2}&\cdots&p^{f_1+\cdots+f_{n-1}}&p^{f_1+\cdots+f_n}\cr
\multispan7\hrulefill.\cr
}}
$$

{\it Remark~4.2}.  More examples of such computations can be found
in \citer\further(Section~4) and a comprehensive summary
in \citer\final().

\medbreak

{\bf 5. Summary}\pointir Let $p$ be a prime number, and let $G$ be a
$p$-group in which every subgroup is the intersection of the index-$p$
subgroups of $G$ containing it.  We show that $G$ has to be
commutative of exponent~$p$.  Let $K$ be a local field with perfect
residue field of characteristic~$p$, and let $L|K$ be a totally
ramified $G$-extension.  If $d_1<d_2<\cdots<d_n$ is the set of
exponents of the differents of the various degree-$p$ extensions of
$K$ in $L$, then the upper ramification breaks of $L|K$ occur
precisely at the integers $t_1<t_2<\cdots<t_n$, where
$d_i=(1+t_i)(p-1)$, and the lower ramification breakes precisely at
$l_1=t_1$ and at $l_i=l_{i-1}+(t_i-t_{i-1})p^{f_1+f_2+\cdots+f_{i-1}}$
($i\in[2,n]$), where $f_j$ is the codimension of $G^{t_j+}$ in the
$\F_p$-space $G^{t_j}$ for every $j\in[1,n]$.

\medbreak {\it It is a pleasure to thank Kirti Joshi for a careful
reading of this Note and for his suggestions for improving the
exposition.  Other readers, if any, owe him their gratitude too.}

\bigbreak
\unvbox\bibbox 

\bye